\documentclass[12pt,a4paper]{article}
\usepackage{amsmath,amsfonts,amssymb,amsthm,amscd}
\usepackage{hyperref}
\usepackage{braket}
\usepackage{epstopdf}

\usepackage{braket}
\usepackage{xypic}
\usepackage{amsmath}

\usepackage[american]{babel}
\usepackage{color} 
\usepackage{epsfig}
\usepackage{caption}
\usepackage{rotating}
\usepackage{setspace}
\usepackage{fancyhdr}
\usepackage{booktabs}
\usepackage{mathrsfs}
\usepackage{amsthm}
\usepackage{amsmath,amssymb}
\usepackage{enumerate}
\usepackage{tikz}
\usepackage{youngtab}
\usepackage{hyperref}
\usepackage{tensor}
\usepackage{mathabx}
\usepackage{booktabs}
\usepackage{mathrsfs}
\usepackage{subfigure}
\usepackage{amsmath,amssymb}

\renewcommand{\[}{\begin{equation}}
\renewcommand{\]}{\end{equation}}

\newtheorem{thm}{Theorem}[section]
\newtheorem{cor}[thm]{Corollary}

\theoremstyle{definition}

\newtheorem{exmp}{Example}[section]

\theoremstyle{remark}

\DeclareMathOperator{\SL}{SL}

\title{On inner-amenability and boundary actions}

\author{J. Bassi }

\begin{document}

\date{}
\maketitle

\begin{abstract}
\noindent Let $\Gamma$ be a discrete countable group. One result in this work is that if $\Gamma$ is ICC inner-amenable non-amenable then it cannot satisfy the (AO)-property, answering a question posed by C. Anantharaman-Delaroche. A generalization of this phenomenon is also considered. It is also proved that if $\Gamma$ is a "sufficiently large" discrete subgroup of a product of locally compact second countable bi-exact groups, then it cannot be inner-amenable. Both these results generalize the well-known fact that ICC non-amenable inner-amenable discrete countable groups cannot be bi-exact.





\end{abstract}

\section{Introduction}
 Topological amenability of boundary actions is a crucial concept in the analytic theory of discrete groups, which proved to capture fine properties of the underlying groups. Given a discrete group $\Gamma$, exactness of $\Gamma$ is equivalent to topological amenability of its left action on its Stone-{\v C}ech boundary $\partial_\beta \Gamma$ (\cite{Oz}). Boundary amenability of certain actions is useful in order to understand the ideal structure of $C^*$-algebras associated to quasi-regular representations (\cite{BeKa, BaRa2}) and is also related to a conjecture by Ozawa on tight inclusions of exact $C^*$-algebras in nuclear $C^*$-algebras (\cite{Oz2, KaKe, BaRa3}). Boundary amenability of the left-right action of $\Gamma \times \Gamma$ on $\partial_\beta \Gamma$, usually referred to as bi-exactness or property $\mathcal{S}$, reflects some form of hyperbolicity and is closely related to the (AO)-property, i.e. temperedness of the left-right representation of $\Gamma \times \Gamma$ in the Calkin algebra of $l^2 (\Gamma)$, first investigated in (\cite{AkOs}). This latter property ensures solidity of the group von Neumann algebra, a remarkable property which implies primeness (\cite{Oz3}). 
 
 In this work we employ regularity properties of certain boundary actions in order to find conditions which deny inner-amenability of a given discrete second countable group. More precisely, using the characterization of the (AO)-property given by the author in \cite{Ba}, it is proven that inner-amenable countable ICC non-amenable groups cannot satisfy the (AO)-property, answering a question posed by C. Anantharaman-Delaroche in \cite{AnDe}. The techniques involved actually show that weaker forms of inner-amenability are denied in case the Calkin representation satisfies weaker forms of (AO). In the second part it is proved, using the notion of bi-exactness for locally compact second countable (lcsc) groups introduced in \cite{De}, that inner-amenable discrete ICC groups cannot be realized as discrete subgroups of products of lcsc bi-exact groups. Both these results generalize the well known fact that discrete countable inner-amenable groups cannot be bi-exact.
 
 \section{Inner-amenability and (AO)-property}
 Let $\Gamma$ be a discrete countable group acting on a countable set $X$. In \cite{Ba} it is introduced the notion of non-standard extension of a boundary in order to translate regularity properties of the associated Calkin representation in the language of measurable dynamics. Non-standard extensions of the Stone-{\v C}ech boundary of $X$ are constructed in the following way: given $\omega \in \partial_\beta \mathbb{N}$, we construct the ultraproduct $C^*$-algebra $\prod_\omega l^\infty X$. Since there is a $\Gamma$-equivariant embeding $l^\infty X \subset \prod_\omega l^\infty X$, by Gelfand duality, there is a surjective $\Gamma$-equivariant continuous map $\phi_\omega : \sigma (\prod_\omega l^\infty X) \rightarrow \Delta_\beta X$. We define the non-standard extension of the Stone-{\v C}ech boundary of $X$ associated to $\omega$ as $\partial_{\beta, \omega} X := \phi_{\omega}^{-1} (\partial_\beta X)$. In \cite{Ba} Definition 2.4 it is introduced a particular set of quasi-invariant measures on $\partial_{ \beta, \omega} X$, which is denoted by $Q-\tilde{\mathcal{P}}^\infty (X)_\omega$. It is proved in \cite{Ba} Theorem 3.2 that the associated $*$-homomorphism $C^*\Gamma \rightarrow \mathbb{B}(l^2 X) / \mathbb{K}(l^2 X)$ factors through the canonical surjection $C^* \Gamma \rightarrow C^*_\lambda \Gamma$ (i.e. the associated group representation is tempered) if and only if every quasi-invariant measure in $Q-\tilde{\mathcal{P}}^\infty (X)_\omega$ gives rise to a tempered Koopman representation.
 
 Let $\Gamma$ be a discrete group and let $C^*_\tau (\Gamma)$ be a $C^*$-completion of the group algebra $\mathbb{C}[\Gamma]$. We say that a unitary representation of $\Gamma$ is $\tau$-continuous if the associated homomorphism $C^* \Gamma \rightarrow \mathbb{B}(H)$ factors through the canonical surjection $C^* \Gamma \rightarrow C^*_\tau \Gamma$. We say that $\Gamma$ is a $\tau$-group if the trivial representation is $\tau$-continuous.

 \begin{thm}
 \label{thm1}
 Let $\Gamma$ be a discrete countable group acting on a countable set $X$ and suppose that the Calkin representation $C^*\Gamma \rightarrow \mathbb{B}(l^2 X)/\mathbb{K}(l^2 X)$ is $\tau$-continuous for a certain $C^*$-topology on $\Gamma$. Let $\Lambda$ be a subgroup of $\Gamma$ which fixes a probability measure on $\partial_\beta X$. Then $\Lambda$ is a $\tau$-group.
 \end{thm}
 \proof By a standard approximation argument, if there is $\Lambda$-invariant probability measure on $\partial_\beta X$, then there is a sequence of almost $\Lambda$-invariant norm-one vectors $(\xi_n)_{n \in \mathbb{N}}$ in $l^2 (X)$ going to zero weakly. For every $n \in \mathbb{N}$ let $\tilde{\xi}_n = \sum_{i \in \mathbb{N}} \alpha_i \xi_n (\gamma_i (\cdot))$ for a choice of $(\alpha_i)_{i \in \mathbb{N}} \subset \mathbb{R}_{+}$ ensuiring that $(\tilde{\xi}_n)$ have uniformly bounded norm in $l^2(X)$. Define $\hat{\xi}_n = \tilde{\xi}_n /\| \tilde{\xi}_n\|_2$. Let $\omega$ be a free ultrafilter on $\mathbb{N}$ and let $\mu$ be the associated quasi-invariant probability measure on $\partial_{\beta, \omega} X$ induced by the sequence $\hat{\xi}_n$. Then the state $\lim_{n \rightarrow \omega}< (\cdot) \xi_n, \xi_n>$ is a vector state coming from the Koopman representation on $L^2(\mu)$ (as shown in \cite{Ba} Proposition 2.13); since this Koopman representation is $\tau$-continuous and this vector is fixed by $\Lambda$, it follows that $\Lambda$ is a $\tau$-group. $\Box$

Recall that a discrete group $\Gamma$ is said to be \textit{inner-amenable} if there is an invariant probability measure on $\Delta_\beta \Gamma \backslash \{e\}$ for the adjoint action of $\Gamma$ given by the extension of $ad (\gamma): x \mapsto \gamma x \gamma^{-1}$, for $x, \gamma \in \Gamma$ (\cite{Ef}). Also, $\Gamma$ satisfies the \textit{(AO)-property} if the left-right representation on its Calkin algebra $\lambda \times \rho : C^*(\Gamma \times \Gamma) \rightarrow \mathbb{B}(l^2 \Gamma)/\mathbb{K}(l^2 \Gamma)$ is tempered. The following question appears in \cite{AnDe}, but the answer is probably already known to experts.
\begin{cor}
Let $\Gamma$ be a discrete countable ICC inner-amenable non-amenable group. Then $\Gamma$ cannot satisfy the (AO)-property.
\end{cor}
\proof This is an application of Theorem \ref{thm1} in the case we consider the left-right action of $\Gamma \times \Gamma$ on $\Gamma$ and the subgroup is the diagonal copy of $\Gamma \subset \Gamma \times \Gamma$. The only thing to observe is that since $\Gamma$ is ICC, the $\Gamma$-invariant probability measure on $\Delta_\beta \Gamma$ is actually supported on $\partial_\beta \Gamma$. $\Box$

\begin{cor}
\label{corsl}
Let $\Lambda$ be a subgroup of $\SL(3,\mathbb{Z})$ fixing a probability measure for the adjoint action on $\partial_\beta \SL(3,\mathbb{Z})$. Then $\Lambda$ has the Haagerup property.
\end{cor}
\proof In \cite{BaRa} it is proved that the $*$-homomorphism $C^* (\SL(3,\mathbb{Z}) \times \SL(3,\mathbb{Z})) \rightarrow \mathcal{Q}(l^2 (\SL(3,\mathbb{Z}))$ induced by the left-right action $(\lambda \times \rho)$ factors through the completion of the group algebra associated to the ideal $c'_0 (\SL(3,\mathbb{Z}) \times \SL(3,\mathbb{Z}))= \{ f \in l^\infty (\SL(3,\mathbb{Z}) \times \SL(3,\mathbb{Z}))\; | \; f(\gamma_n, \eta_n) \rightarrow 0 \; \mbox{ for } (\gamma_n, \eta_n) \rightarrow (\infty, \infty)\}$. The result follows since the restriction of this ideal to the diagonal copy of $\SL(3,\mathbb{Z})$ in $\SL(3,\mathbb{Z}) \times \SL(3,\mathbb{Z})$ coincides with $c_0 (\SL(3,\mathbb{Z}))$. $\Box$

It is an open problem whether $\SL(3,\mathbb{Z})$ satisfies the (AO)-property. If it does, then the subgroup $\Lambda$ in Corollary \ref{corsl} would actually be forced to be amenable.

In \cite{BeKa} the authors introduced different classes of subgroups of a given discrete group. In particular, if $\Lambda$ is a subgroup of a discrete group $\Gamma$, then $\Lambda$ is said to have the \textit{spectral gap property} (as a subgroup of $\Gamma$) if the associated representation of $\Lambda$ on $l^2 (\Gamma / \Lambda \backslash \{e\Lambda\})$ does not weakly contain the trivial representation. Note that the above definition is equivalent to the requirement that there is no $\Lambda$-invariant probability measure on $l^\infty (\Gamma / \Lambda \backslash \{ e \Lambda\})$.

\begin{cor}
\label{cor1}
Let $\Lambda$ be a non-amenable subgroup of a discrete countable group $\Gamma$. If the Calkin representation $\Lambda \rightarrow \mathbb{B}(l^2 (\Gamma /\Lambda))$ is tempered, then $\Lambda$ has the spectral gap property if and only if every orbit of the form $\Lambda (g\Lambda)$, with $g \notin \Lambda$, is infinite.
\end{cor}
\proof By Theorem \ref{thm1} there is no $\Lambda$-invariant probability measure on $\partial_\beta (\Gamma / \Lambda)$. Hence the only $\Lambda$-invariant probability measures on $l^\infty (\Gamma /\Lambda \backslash \{e\Lambda\})$ can arise as probability measures on $\Gamma /\Lambda \backslash \{e \Lambda\}$. Such a measure exists if and only if there is a finite $\Lambda$-orbit on $\Gamma /\Lambda \backslash \{e\Lambda\}$. $\Box$

An example of such a situation is the inclusion of $\SL(2,\mathbb{Z})$ inside $\SL(3,\mathbb{Z})$ given by
\begin{equation*}
\left(\begin{array}{cc} 1 & 0 \\ 0 & \SL(2,\mathbb{Z}) \end{array}\right) \subset \SL(3,\mathbb{Z}).
\end{equation*}
By the way the conclusion that this subgroup has spectral gap (which was already observed in \cite{BeKa}), actually follows directly from the fact proved in \cite{BaRa3} that the action of $\SL(3,\mathbb{Z})$ on the Stone-{\v C}ech boundary of $\SL(3,\mathbb{Z})/\SL(2,\mathbb{Z})$ is topologically amenable.

\section{Discrete subgroups of products of bi-exact groups}
We begin this section by fixing some notation and refer to \cite{De} for a detailed survey of bi-exactness for lcsc groups. Let $G$ be a lcsc group. There is a universal small at infinity left-right-equivariant compactification of $G$ given by the spectrum of the $C^*$-algebra $C(h^u G) = \{ f \in C_b^u (G) \; | \; \rho_g f - f \in C_0 (G) \; \mbox{ for all } g \in G\}$, where $C_b^u (G)$ is the algebra of bounded functions on $G$ such that $\|\lambda_g f - f\|_\infty \rightarrow 0$ and $\| \rho_g f-f\|_\infty \rightarrow 0$ whenever $g \rightarrow e$, for $f \in C_b^u (G)$. The boundary of $h^u G$ is denoted by $\nu^u G$. Both the left and the right action of $G$ on $G$ extend to $\nu^uG$ and the right action extends to the trivial action (i.e. the compactification is small at infinity). Recall that $G$ is bi-exact if its left action on $\nu^u G$ is topologically amenable.

In \cite{HaSk} and \cite{ChSiBu} the authors prove that certain lattices in products of topological groups fail to be inner amenable. The following is a related result.
\begin{thm}
Let $\Gamma$ be a discrete countable ICC subgroup of a product of locally compact second countable non-amenable bi-exact groups $G_1 \times G_2 \times ... \times G_n$ such that the projection of $\Gamma$ in each factor contains a non-amenable discrete subgroup of $G_i$, or that $G$ is irreducible, i.e. the projection on each factor is dense. Then $\Gamma$ is not inner-amenable.
\end{thm}
\proof For every $i \in \{1,..., n\}$ let $\phi_i : \Delta_\beta \Gamma \rightarrow h^u G_i$ be the continuous $\Gamma \times \Gamma$-equivariant map extending the projection on the $i$-th factor (the action of $\Gamma \times \Gamma$ on each $h^u G_i$ is given by $(\gamma, \eta, x_i) \mapsto \phi_i (\gamma) x_i \phi_i (\eta)^{-1}$). Clearly $\phi_i^{-1} \nu^u G_i \subset \partial_\beta \Gamma$ for every $1 \leq i \leq n$. Moreover, $\bigcup_{i=1,...,n} \phi_i^{-1} \nu^u G_i = \partial_\beta \Gamma$, indeed, suppose this is not the case, then we can find a point $x=(x_1,...,x_n) \in G_1 \times ... \times G_n$ and $\omega \in \partial_\beta \Gamma$ such that $\phi_i (\omega)=x_i$ for every $1 \leq i \leq n$. But then there is a compact set $K$ with non-empty interior around $x$ such that $A_K := \{ \gamma \in \Gamma \; | \; \gamma \in K\} \in \omega$, but since $K$ is compact and $\Gamma$ is discrete, $A_K$ is finite, contradicting the fact that $\omega \in \partial_\beta \Gamma$. Now we want to show that for every $1 \leq i \leq n$, $\omega \in \phi_i^{-1} \nu^u G_i$ and $\gamma \in \Gamma$ we have $\phi_i (\gamma \omega \gamma^{-1})= \gamma \phi_i (\omega)$. First observe that each compact set $\phi_i^{-1} \nu^u G_i$ is $\Gamma \times \Gamma$-invariant in virtue of the $\Gamma \times \Gamma$-equivariance of $\phi_i$. Hence for every $i$ we have $\phi_i (\gamma \omega \gamma^{-1})=\gamma \phi_i (\omega) \gamma^{-1} = \gamma \phi_i (\omega)$ since the right action on $G_i$ extends trivially to $\nu^u G_i$. Suppose now that $\mu$ is an $ad(\Gamma)$-invariant probability measure on $\partial_\beta \Gamma$. Then there is at least one $i$ for which $\mu_i =\mu|_{\phi_i^{-1} \nu^u G_i } >0$. Up to normalization, for this choice of $i$, $\mu_i$ is an $ad(\Gamma)$-invariant probability measure on $\phi_i^{-1} \nu^u G_i $. We can push this to a probability measure on $\nu^u G_i$, which, by the above discussion, is fixed by the left action of $\Gamma$. Now, if $\Gamma$ is not irreducible, we employ the requirement that $\phi_i (\Gamma)$ contains a non-amenable discrete subgroup of $G_i$ in order to reach a contradiction. If $\Gamma$ is irreducible, it follows from the continuity of the action of $G_i$ on $\mathcal{P} (\nu^u G_i)$ that $G_i$ fixes a probability measure on $\nu^u G_i$; now, since $G_i$ acts amenably on $\nu^u G_i$, it follows that it acts amenably on $\mathcal{P}(\nu^u G_i)$, contradicting the existence of a fixed probability measure. $\Box$

\begin{exmp}
For every prime number $p$, the group $\SL(2,\mathbb{Z}[1/p])$ is an irreducible lattice in $\SL(2,\mathbb{R}) \times \SL(2,\mathbb{Q}_p)$, hence it is not inner-amenable.
\end{exmp}

\section{Acknowledgments}

The author aknowledges the support of the grant SOE-Young Researchers 2024 ”Generalized Akemann-Ostrand Property: Analytical, Dynamical And Rigidity Properties (GAOPADRP)”, CUP: E83C24002550001 and the MIUR Excellence Department Project 2023–2027 MatMod@Tov awarded to the Department of Mathematics, University of Rome Tor Vergata. He also akcnowledges the support of INdAM-GNAMPA and of IM PAN, where this work was partially developed. He thanks  Prof. R. Conti , Prof. K. Li and Prof. F. R\u{a}dulescu for interesting discussions on the topic of this work.


\baselineskip0pt
\bigskip
  \footnotesize

  J. Bassi, \textsc{Department of Mathematics, University of Tor Vergata, Via della Ricerca Scientifica 1, 00133, Roma, Italy}\par\nopagebreak

\end{document}